\documentclass[11pt]{article}
\usepackage{times,amsmath,amsthm,amssymb,graphicx,hyperref,float,xcolor}
\usepackage{fancyhdr}     
\usepackage{setspace}    
\usepackage{booktabs}    
\usepackage{longtable}   
\usepackage{multirow}
\usepackage{textcomp}    
\usepackage{bm}         
\usepackage{pdfpages}    
\newtheorem{defn}{Definition}
\newtheorem{rem}{Remark}
\newtheorem*{pf}{Proof}
\newtheorem{lem}{Lemma}
\newtheorem{thm}{Theorem}
\newtheorem{exmp}{Example}

\newcommand{\disp}{\displaystyle}

\newcommand{\bx}{\mathbf{x}}
\newcommand{\by}{\mathbf{y}}
\newcommand{\bbeta}{\bm{\beta}}
\newcommand{\bxi}{\bm{\xi}}
\begin{document}
\global\def\refname{{\normalsize \it References:}}
\baselineskip 12.5pt
%
%
%
\title{\LARGE \bf Identification of the Source for Full Parabolic Equations}

\date{}

\author{\hspace*{-10pt}
\begin{minipage}[t]{2.3in} \normalsize \baselineskip 12.5pt
\centerline{GUILLERMO FEDERICO UMBRICHT}
\centerline{Instituto de Ciencias e Instituto del Desarrollo Humano, Univ.~Nac. de Gral.~Sarmiento}
\centerline{J. M. Guti\'errez 1150, Los Polvorines, Buenos Aires, ARGENTINA}
\centerline{Centro de Matem\'atica Aplicada, Escuela de Ciencia y Tegnolog\'ia, Univ.~Nac. de Gral.~San Mart\'in}
\centerline{25 de Mayo y Francia, San Mart\'in, Buenos Aires, ARGENTINA}
\end{minipage}
%
\\ \\ \hspace*{-10pt}
\begin{minipage}[b]{6.9in} \normalsize
\baselineskip 12.5pt {\it Abstract:}
In this work, we consider the problem of identifying the  time independent source for full parabolic equations in $\mathbb{R}^n$ from noisy data. This is an ill-posed problem in the sense of Hadamard. 
To compensate the factor that causes the instability, a family of parametric regularization operators is introduced, where the rule  to select the value of the regularization parameter is included. This rule, known as regularization parameter
choice rule, depends on the data noise level and the degree of smoothness that it is assumed for the source. The proof for the stability and convergence of the regularization criteria is presented and a H\"older type bound is obtained for the
estimation error. Numerical examples are included to illustrate the effectiveness of this regularization approach. 
\\ [4mm] {\it Key--Words:}
inverse and ill-posed problem, regularization operator, transport equation, Fourier transform.
\end{minipage}
\vspace{-10pt}}

\maketitle

\thispagestyle{empty} \pagestyle{empty}
%
%
\section{Introduction}
\label{s:1} \vspace{-4pt}

The problem of source identification has been studied and analyzed in different areas of applied mathematics for the last decades. It has also received considerable attention from many current researchers in science and engineering. Applications can be found in problems related to heat conduction \cite{Elden00,ZhaoMeng2011},fissure identification \cite{Zeng96}, geophysical prospecting \cite{BEROZA88}, contaminant detection \cite{Letall2006} and tumor cell detection \cite{Macleod99}, to mention some.

The determination of sources turns out to be an ill-posed problem in the sense of Hadamard \cite{Hadamard23} since the solution does not depend continuously on the data. 

Among the most significant tools used to determine a source, one can find in the literature the potential logarithmic method \cite{Ohe94}, the projective method \cite{Nara03}, the Green function \cite{Hon10}, dual reciprocity boundary element methods \cite{FARCAS03}, the dual reciprocity method \cite{Sun97}, the fundamental MFS solution method \cite{Jin07} and method by using the curve L \cite{HOL93}.

Regarding the transport term in parabolic differential equations, there are no many articles published for the general case, most of the papers available in the literature focused on the heat equation. Sources of the heat equation are recovered using different methods and strategies, see for instance \cite{FARCAS06,AHMADABADI09,JM07}. There are many articles that analyze particular cases with simplifications or restrictions on the mathematical models, the type of source, the border conditions or the chosen domain, like  \cite{FARCAS06, Liu09, Savateev95, Cannon98, Yaetal2008, AHMADABADI09, JM07}. The most commonly used methods are limit element method \cite{FARCAS06, JM07}, MFS fundamental solution method \cite{Yaetal2008, AHMADABADI09}, Ritz-Galerkin method \cite{Ras11}, differences method finite \cite{Yaetal2010}, non-mesh method \cite{Yaetal2009}, conditional stability method \cite{Yama93} and the firing method \cite{Liu09}. 

On the other hand, regularization methods \cite{Engel1996, Kirsch2011, Mazzieri14} play an important role in the estimation of unstable solutions. The most widely used approaches are the iterative regularization method \cite{JM08}, the simplified method of regularization Tikhonov \cite{Fu04, YangFu2010, Cheng07, Cheng08}, the modified regularization method \cite{YangFu2011, YF10, ZhaoMeng2014}, Fourier truncation \cite{YangFu2011}, the method of mollification \cite{YangFu2014}. In particular for the determination of the source for a parabolic equation, certain regularization techniques are applied to specific transport equation. In \cite{Sivergina2003} the authors focused on a convection-diffusion equation, while in \cite{DF2009, DFY2009, Elden00, TLA2005, TQA2006, YF10, Yaetal2010, ZhaoMeng2014} only diffusion is considered. 

Recently, the problem to find the source term as in this manuscript was solved by the quasi-reversibility method, see \cite{Ngu2019}. This method can be used to solve inverse source problem for nonlinear parabolic equations, the difference of this manuscript and that paper is the observation data \cite{Le2020}.

This work aims to the determination, from noisy measurements taken at an arbitrary fixed time, of the real-valued function of n real variables, independent of time, in an evolutionary equation of transport in an unbounded domain. 
This is an ill-posed problem 
because the high frequency components of arbitrarily small data errors can lead to arbitrarily large errors in the solution \cite{Engel1996, Kirsch2011}.

Here, a family of regularization operators is designed to compensate the factor that causes the instability of the inverse operator. The parametric regularization operators lead to a family of well-posed problems that approximates the given ill-posed problem. The regularization operator family proposed here turns out to be an $n$-dimensional generalization of the modified regularization method considered in  \cite{QianFu06,QianFu07,Xiao11,YangFu2012,ZhaoMeng2014}. In these articles the authors estimate the source of the one-dimensional equation of heat from data measured in a fixed moment of time $(t = 1)$ by adding a penalizing term and the parameter choice rule depends on the norm of the unknown function. In contrast, this paper analyzes the general n-dimensional parabolic equation while relaxing the conditions on the assumptions for the regularization process.\\

The stability and convergence of this regularization method are analyzed and a H\"older type bound  is derived for the estimation error. In order to  illustrate the regularization performance, some numerical examples for the 1D, 2D and 3D cases are included. 
\section{The ill-posed mathematical framework}\label{s:2} \vspace{-4pt}

We consider the  problem of determining the source $f$ for the following parabolic equation 
\begin{equation}
\label{transpeqn}
\begin{cases} 
u_t({\bx},t)=\alpha^2 \Delta u(\bx,t)-\bbeta \cdot \nabla (u(\bx,t))-\nu u(\bx,t)+f(\bx), \qquad & \bx \in \mathbb{R}^{n}, \, t>0, \\
u(\bx,0)=0, \qquad & \bx \in \mathbb{R}^{n}, \\
u(\bx,t_0)=y(\bx), \qquad & \bx \in \mathbb{R}^{n}, \, t_0>0, \\ 
\end{cases}
\end{equation}
where $\alpha^2, \nu >0$, \, $\bbeta\in \mathbb{R}^{n}$ are given, $ \Delta $ denotes the Laplacian operator, $\nabla$ denotes the Nabla operator and $``\cdot"$ is the usual inner product in $\mathbb{R}^{n}$. Note that this is a linear parabolic equations with constant coefficients. The existence and uniqueness of the solution to \eqref{transpeqn} is discussed in \cite{Pao1998}. 
It is assumed that $u(\cdot, t), f(\cdot) \in L^{2}(\mathbb{R}^{n})$ are unknown functions and that $y \in L^2(\mathbb{R}^{n})$ can be measured with certain noise level $\delta$, i.e., the data function $y_{\delta} \in L^{2}(\mathbb{R}^{n})$ satisfies 
\begin{equation}
\label{noiselevel}
    ||y-y_{\delta }||_{L^{2}(\mathbb{R}^{n})}\le \delta, \qquad 0<\delta\le \delta_M,\\
\end{equation} 
where $\delta_M \in \mathbb{R}_{>0} $ represents the maximum level of noise. In practice , $\delta_M$ may be estimated from the error committed by the measuring instruments.

The analysis of the equation with boundary values and initial conditions in \eqref{transpeqn} is perform in the frequency space. 

\begin{defn} 
Let $g \in L^{2}(\mathbb{R}^{n})$. The Fourier n-dimensional transform is defined 
\begin{equation}
\label{deftransf}
\hat{g}(\bxi) :=\left(\frac{1}{\sqrt{2\pi}}\right)^{n} \int\limits_{\mathbb{R}^n}{e^{-i \bxi \cdot \bx} \, g (\bx) \,  d\bx}, \qquad  \qquad  \bxi, \bx \in \mathbb{R}^n.
\end{equation}
\end{defn}

By using the above definition \eqref{deftransf}, the system \eqref{transpeqn} can be written in the frequency space as
\begin{equation}
\label{eqnfreq}
\begin{cases}
\hat{u}_{t}(\bxi ,t)=- z(\bxi) \, \hat{u}(\bxi,t)+\hat{f}(\bxi),\qquad  & \bxi \in \mathbb{R}^{n}, \, t>0, \\
 \hat{u}(\bxi ,0)=0, & \bxi \in \mathbb{R}^{n}, \\
\hat{u}(\bxi ,t_0)=\hat{y}(\bxi ), & \bxi \in \mathbb{R}^{n}, \, t_0>0,
\end{cases}
\end{equation}
where 
$$z(\bxi)=\alpha^2 \| \bxi\|^2 + i\,\bbeta \cdot \bxi + \nu \, \in \mathbb{C}.$$
Solving \eqref{eqnfreq} in the frequency space, we obtain the solution  
\begin{equation}
\label{solutionu}
\hat{u}(\bxi,t) = \frac{1-e^{-z(\bxi)\, t}}{z(\bxi)}\hat{f}(\bxi ).
\end{equation}
Since $\hat{u}(\bxi ,t_0)=\hat{y}(\bxi )$, an expression for the source  in the frequency space  is obtained by evaluating the equation \eqref{solutionu} in $t=t_0$, that is,   
\begin{equation}
\label{illf}
\hat{f}(\bxi )=\Lambda (\bxi )\hat{y}(\bxi ),
\end{equation}
 where 
\begin{equation}
\label{Lambda}
 \Lambda (\bxi )=\frac{z(\bxi)}{1-e^{-z(\bxi)\, t_0}}.
\end{equation}
Denoting  $\hat{f_{\delta}}(\bxi )=\Lambda (\bxi )\hat{y_{\delta}}(\bxi )$, we have \\
$$ \|\hat{f} - \hat{f_{\delta}}\|_{L^2(\mathbb{R}^n)} = \| \Lambda(\bxi) ( \hat{y}(\bxi )- \hat{y_{\delta}}(\bxi ))\|_{L^2(\mathbb{R}^n)}=
\| \Lambda(\bxi)\|_{L^2(\mathbb{R}^n)} \| \hat{y}(\bxi)- \hat{y_{\delta}}(\bxi)\|_{L^2(\mathbb{R}^n)}.$$

Since
\begin{eqnarray}
  \label{lim}
 | \Lambda (\bxi )|
&=&\left| \frac{  z(\bxi)}{1-e^{-(\alpha^2 \| \bxi\|^2 + i\,\bbeta \cdot \bxi + \nu)\, t_0}}\right| \geq 
 \frac{  |z(\bxi)|}{1+e^{-(\alpha^2 \| \bxi\|^2 + \nu)\, t_0} \left| e^{-i\,\bbeta \cdot \bxi \, t_0}\right| } \nonumber \\
&\geq& 
  \frac{  |\alpha^2 \| \bxi\|^2+ \nu + i\,\bbeta \cdot \bxi  |}{1+e^{-(\alpha^2 \| \bxi\|^2 + \nu)\, t_0} },
\end{eqnarray}
$\Lambda (\bxi)$ increases without bound as $\|\bxi\| \to \infty$  amplifying the high frequency components of the observation error $ \hat{y}(\bxi) - \hat{y}_{\delta} (\bxi)$. This fact might lead to a large estimation error $ \|\hat{f} - \hat{f_{\delta}}\|_{L^2(\mathbb{R}^n)}$ even for small observation errors, hence one of the Hadamard conditions is not satisfied \cite{Hadamard23}.

\section{Regularization operators}\label{s:3}\vspace{-4pt}

In this section we propose a regularization operator taking into account the inestability factor in the inverse operator.
We notice that the resulting operator is equivalent to the one that is obtained by using the quasi-reversibility method \cite{Lat1969}.
Basic theoretical issues related to regularization operators are included, more information can be found  in \cite{Engel1996, Kirsch2011}.

\begin{defn} 
Let $T : Y \longrightarrow X$,  $X$ and $Y$ be Hilbert spaces and $T$ be an unbounded operator. A regularization strategy for $T$ is a family of linear and bounded operators 
\begin{equation}
\label{DefinitionR}
R_{\mu} : Y \longrightarrow X, \quad \mu>0, \quad /  \quad  \lim _{\mu \to 0^+} R_{\mu} y = Ty, \quad \forall y \in Y.\\
\end{equation}
\end{defn}

Let us define the parametric family of linear operators 
$R_\mu: L^{2}(\mathbb{R}^{n}) \to L^{2}(\mathbb{R}^{n})$ for $\mu \in \mathbb{R}_{>0}$, by 
\begin{equation}
\label{familyR}
R_{\mu} \hat{y}(\bxi) :=  \frac{\Lambda (\bxi)}{1+\mu^2 \|\bxi\|^2} \, \hat{y}(\bxi),
\end{equation}
where $\Lambda (\bxi )$ is given in \eqref{Lambda} and $\mu$ is the regularization parameter. 
Note that the denominator in \eqref{familyR} was introduced for stabilization purposes. The properties for the operator family $\{ R_{\mu},\, \mu>0 \}$ are stated in the following theorem.

\begin{thm}
Let us consider the problem of identifying $f$ from noisy data $ y_{\delta}(x)$  measured at a given time instant $t_0>0$, where  $\delta$ is the noise level defined in \eqref{noiselevel}. Let the functions $u$ and $f$ satisfy the following differential equation with initial condition 
\begin{equation}
\label{illpp}
\begin{cases} 
u_t(\bx,t)=\alpha^2 \Delta u(\bx,t)-\bbeta \cdot \nabla (u(\bx,t))-\nu u(\bx,t)+f(\bx), \qquad & \bx \in \mathbb{R}^{n}, \, t>0, \\
u(\bx,0)=0, \qquad & \bx \in \mathbb{R}^{n}, \\
\end{cases}
\end{equation}
and let $\{R_{\mu}\}$ be the family of operators defined in \eqref{familyR}.
Then, for every $y(\bx)=u(\bx,t_0)$ there exists an a-priori parameter choice rule for $\mu >0$ such that the pair $(R_{\mu}, \mu)$ is a convergent regularization method for solving the identification problem \eqref{illpp}.
\end{thm}
\begin{pf}
The factor $\disp \frac{\Lambda (\bxi)}{1+\mu^2\|\bxi\|^2}$ is bounded for all $\bxi$  since it is continuous for all $\bxi \in \mathbb{R}^{n}$ and
\begin{equation*}
\begin{split}
\lim _{\|\bxi\| \to  \infty} \left\vert \frac{\Lambda (\bxi)}{1+\mu^2\|\bxi\|^2} \right\vert 
=&\lim _{\|\bxi\| \to  \infty} \left\vert \frac{ \alpha^2 \|\bxi\|^{2} + i \bbeta \cdot \bxi + \nu}{(1-e^{-(\alpha^2 \|\bxi\|^{2} + i \bbeta \cdot \bxi + \nu) \, t_0})(1+\mu^2\|\bxi\|^2)}\right\vert \\
 =&\frac{\alpha^2 }{\mu^2}<\infty. 
\end{split}
\end{equation*}

Hence, for $\mu>0$, $R_{\mu}$ is a linear continuous operator and we have that 
$$\lim _{\mu \to 0^+} R_{\mu} \hat{y} = \Lambda \hat{y},$$
for  $ \hat{y} \in L^{2}(\mathbb{R}^{n})$, then $R_{\mu}$ is a regularization strategy for $\Lambda$. Therefore, by Proposition 3.4 in \cite{Engel1996}, for $y(\bx)=u(\bx,t_0)$ there exists an a-priori parameter choice rule $\mu$ such that $(R_{\mu}, \mu)$ is a convergent regularization method for solving \eqref{illf}. The regularized solution to the inverse problem in the frequency space is given by
\begin{equation}
\label{ffregtransf}
\hat{f}_{\delta ,\mu} = \frac{\Lambda (\bxi)}{1+\mu^2 \|\bxi\|^2} \, \hat{y}_\delta (\bxi).
\end{equation}
Therefore, an estimated function to $f$ in \eqref{illpp} is given by the expression
\begin{equation}
\label{ffreg}
f_{\delta ,\mu} =\left(\frac{1}{\sqrt{2\pi}}\right)^{n} \int\limits_{\mathbb{R}^n} e^{i \bxi \cdot \bx} R_{\mu} \hat{y}_{\delta} (\bxi) \,  d\bxi  =\left(\frac{1}{\sqrt{2\pi}}\right)^{n} \int\limits_{\mathbb{R}^n}{e^{i \bxi \cdot \bx} \frac{\Lambda (\bxi)}{1+\mu^2 \|\bxi\|^2} \, \hat{y}_\delta (\bxi) \,  d\bxi}.
\end{equation}
\end{pf}

\section{Error analysis}\label{s:4} \vspace{-4pt}

In order to analyze the regularization performance, we first introduce some results that will be used later to obtain a bound for the  error between the source $f(\bx)$ and its estimate $f_{\delta ,\mu }(\bx)$, that we will be referred to as the regularization error. 
\begin{lem}
\label{lemma1}
For $\omega \in \mathbb{C}$ with $Re(\omega )>0$ holds 
$\disp \left\vert \frac{1}{1-e^{-\omega }} \right\vert \le \frac{1}{1-e^{-Re(\omega )}}.$
\end{lem}

\begin{pf}
Euler's formula for a complex number $\omega=a+bi$,  and the parity for sine and cosine yields
 $e^{-(a+bi) }=e^{-a} cos(b) - i e^{-a} sin(b)$. Then, adding and subtracting $2 e^{-a}=2 e^{-Re(\omega)}$, after algebraic operations one gets
\begin{equation*}
\disp{ |  1-e^{-\omega } |^2
				=(1-e^{-Re(\omega )})^2+2e^{-Re(\omega )}(1-cos(Im(\omega)))} 
				\geq (1-e^{-Re(\omega )})^2.
\end{equation*}

Hence
 $$\disp |  1-e^{-\omega } | \geq 1-e^{-Re(\omega )}  \Longrightarrow \disp \left\vert \frac{1}{1-e^{-\omega }} \right\vert \le \frac{1}{1-e^{-Re(\omega )}},$$
and the proof is completed
\end{pf}

\begin{lem} 
\label{lemma3}
The function $f:\mathbb{R}_{>0} \to \mathbb{R}$ given by 
$\disp f(x)=
\begin{cases}
  \disp \frac{x}{1-e^{-x}}, & 0<x<1,\\
  \disp \frac{1}{1-e^{-x}}, & 1\leq x,
  \end{cases}$
satisfies $f(x) \leq 2 $.
\end{lem}
\begin{pf}
First, let us consider the function $f$ in $(0,1)$. Differentiating, in this case, we have
$\disp  f'(x)= \left(\frac{x}{1-e^{-x}}\right) '=\disp \frac{1+e^{-x}(-1+x)}{(1-e^{-x})^2} > 0.$ Then the function $f$ is increasing in $(0,1)$ and $\disp \frac{x}{1-e^{-x}} \leq \frac{1}{1-e^{-1}}.$\\
On the other hand, for $x>1$,  $f'(x)=\disp \left (\frac{1}{1-e^{-x}}\right) '=\disp \frac{-e^{-x}}{(1-e^{-x})^2} < 0,$ then  $ f $ is decreasing $ \forall x>1 $ and $\disp \frac{1}{1-e^{-x}} \leq \frac{1}{1-e^{-1}}.$\\
Therefore we have that $\disp f(x)< \frac{1}{1-e^{-1}} \, \, \forall x>0$ and since $\disp \frac{1}{1-e^{-1}} \leq 2$ the proof is completed.
\end{pf}

\begin{lem} 
\label{lemma2}
 Let $\rho \in \mathbb{R}$.  If $0<\mu<1 $ we have $\disp \frac{\vert \rho \vert}{1+\rho^2\mu^2} \le \frac{1}{2\mu}.$ 
 Moreover, for $\, \alpha^{2}, \nu >0$  the following inequality  holds  $\disp \frac{\alpha^2\rho^2+\nu }{1+\rho^{2}\mu^2} 
  \le \max \left\{ \nu, \frac{\alpha^{2} }{\mu^{2} } \right\}.$
\end{lem}

\begin{pf}
Since $a^2+b^2 \geq 2ab$ for all $a,b \in \mathbb{R}$, we have that
$$ 1+ \vert \rho \vert ^2 \mu ^2 \geq 2 \vert \rho \vert \mu \Longrightarrow \disp \frac{\vert \rho \vert}{1+\rho^2\mu^2} \le \frac{1}{2\mu}. $$

Now, let $\disp k(\rho)=\frac{\alpha^2\rho^2+\nu 
}{1+\rho^2\mu^2}$, then 
$\disp k'(\rho )=\frac{2\rho (\alpha^2-\nu \mu^2)}{(1+\rho^2\mu^2)^{2}}$ and $k$ has only one critical point at $\rho =0$. Consider three cases: 
\begin{itemize}
\item $\alpha^2=\nu \mu^2$: we have  $k(\rho )=\nu$, constant $\forall \rho \in \mathbb{R}$.
\item $\alpha^2<\nu \mu^2$: then the function $k$ reaches its global maximum value $\nu$ at $\rho=0$.\\ 
\item $\alpha^2>\nu \mu^2$: since $k$ is an even function  and it is increasing for $\rho>$0 with 
 $\displaystyle \lim_{\rho \to \pm \infty} k(\rho)=\frac{\alpha^{2}}{\mu^{2}}$ we have $k(\rho)\leq \frac{\alpha^{2}}{\mu^{2}}$,
 \end{itemize}
 and the proof is completed. 
 \end{pf}

\begin{lem} 
\label{lemma4}
For $\alpha^2,\nu,t_0 >0$, \, $\bbeta\in \mathbb{R}^{n}$ 
 and $0<\mu <1 $  we have that
\begin{center}
$ \disp  \left\vert  \frac{\Lambda(\bxi) } {1+\mu^2\| \bxi\|^2 }  \right\vert \le 
\frac{2}{\mu^2} M $. 
\end{center}
where $M=\max \left\{\frac{1}{t_0} + \frac{\sqrt{n} \|\bbeta\|_\infty}{2\nu t_0} ;\nu 
+\alpha^2+\frac{\sqrt{n} \|\bbeta\|_\infty}{2} \right\}$ and $\|\bbeta\|_\infty=\max\limits_{1\leq j \leq n} |\bbeta_{j}|.$
\end{lem} 
\begin{pf}
From equation \eqref{Lambda} and Lemma \ref{lemma1} we have   
\begin{equation}
\label{D_ineq}
	\left\vert \frac{\Lambda (\bxi)}{1+\mu^2 \|\bxi\|^2}\right\vert \le 
\frac{\alpha^2 \|\bxi\|^2 + i \bbeta \cdot \bxi + \nu}{(1-e^{-(\alpha^2 \|\bxi\|^2 + \nu) \, t_0})(1+\mu^2\|\bxi\|^2) }.
\end{equation}
\setlength{\leftskip}{0pt} 
\setlength{\leftskip}{0pt}
 
\hspace{-1.2cm} 
\emph{If $(\alpha^2\|\bxi\|^2 + \nu) \, t_0 \ge 1$ } : Using the triangular inequality, Lemmas \ref{lemma3} and \ref{lemma2} 
\setlength{\leftskip}{1cm}

\begin{equation} \label{part1}
\begin{split}
\phantom{space} \frac{\alpha^2 \|\bxi\|^2 + \nu +  \vert \bbeta \cdot \bxi \vert}{(1-e^{-(\alpha^2\|\bxi\|^2+\nu ) \, t_0})(1+\|\bxi\|^2\mu^2)} 
&\leq 2 \left( \frac{\alpha^2 \|\bxi\|^2 + \nu }{1 + \|\bxi\|^2 \mu^2} + \frac{\sqrt {n} \|\bbeta\|_\infty \vert \|\bxi\| \vert }{1 + \|\bxi\|^2 \mu^2} \right) \\
&\leq 2 \max \left\{\nu, \frac{\alpha^2}{\mu^2} \right\} +\frac{\sqrt {n} \|\bbeta\|_\infty}{\mu} \\
&\le  \frac{2}{\mu^2}\left( \nu +\alpha^2+\frac{\sqrt {n} \|\bbeta\|_\infty }{2} \right). 
\end{split}
\end{equation}

\setlength{\leftskip}{0pt} 
\setlength{\leftskip}{0pt}

\emph{If $(\alpha^2\|\bxi\|^2+\nu) \, t_0 \in (0,1)$}: Observe that

\setlength{\leftskip}{1cm}

\begin{eqnarray*}
\phantom{space}
\frac{\alpha^2\|\bxi\|^2+\nu + \vert \bbeta \cdot \bxi \vert}{(1-e^{-(\alpha^2\|\bxi\|^2+\nu ) \, t_0})(1+\|\bxi\|^2\mu^2)} 
= 
\frac{(\alpha^2\|\bxi\|^2+\nu) \, t_0}{(1-e^{-(\alpha^2\|\bxi\|^2+\nu ) \, t_0})(1+\|\bxi\|^2\mu^2) \, t_0} \\
\\
\hspace*{2cm} + \frac{\vert  \bbeta \cdot \bxi \vert \, (\alpha^2\|\bxi\|^2+\nu) \, t_0}{(1-e^{-(\alpha^2\|\bxi\|^2+\nu ) \, t_0})(1+\|\bxi\|^2\mu^2)\, (\alpha^2\|\bxi\|^2+\nu) \, t_0}.
\end{eqnarray*}
Using Lemmas \ref{lemma3} and \ref{lemma2}

\begin{equation}
\label{part2}
\begin{split}
\frac{\alpha^2\|\bxi\|^2+\nu + \vert \bbeta \cdot \bxi \vert}{(1-e^{-(\alpha^2\|\bxi\|^2+\nu ) \, t_0}) (1+\|\bxi\|^2 \mu^2)} & \\
\\
 &
\hspace{-2cm} \le 2  \left(\frac{1}{(1+\|\bxi\|^2\mu^2) \, t_0}  +   \frac{\vert \bxi \cdot \bbeta \vert}{(1+\|\bxi\|^2 \mu^2)(\alpha^2 \|\bxi\|^2+\nu) \, t_0} \right)\\
\\
&
\hspace{-2cm}\le  2  \left(\frac{1}{(1+\|\bxi\|^2\mu^2) \, t_0}  +   \frac{\sqrt{n} \|\bbeta\|_\infty }{2 \mu (\alpha^2\|\bxi\|^2+\nu) \, t_0} \right) \qquad \qquad \\
\\
 &
\hspace{-2cm}\le  2  \left(\frac{1}{t_0}  +   \frac{\sqrt{n}\|\bbeta\|_\infty}{2 \mu \nu t_0} \right)
	\le    \frac{2}{\mu^2 t_0}  \left(1  +   \frac{\sqrt{n} \|\bbeta\|_\infty}{2 \nu} \right). 
	\end{split}
	\end{equation}
	\end{pf}

\begin{defn} 
The norm in the Sobolev space $H^{p}(\mathbb{R}^{n}), \, p>0$ is defined as follows
\begin{equation}
\label{Boundf}
\|  f\|  _{H^{p}(\mathbb{R}^{n})} := \left(\,\, \int\limits_{\mathbb{R}^n}{| \hat{f}|^2 \left(1+\|\bxi\|^2 \right)^{p}\, d\bxi} \,\, \right)^{1/2}.
\end{equation}
\end{defn}

\setlength{\leftskip}{0pt}
\begin{thm}
\label{boundestimate1}
Consider the inverse problem of determining the source $f(\bx)$ in (\ref{transpeqn})-(\ref{noiselevel}). Let $f_{\delta ,\mu }(\bx)$ be the regularization solution  given in (\ref{ffreg}) and assume that $\|  f \| _{H^p(\mathbb{R}^n)}<C $ ( $f$ is bounded in $H^p(\mathbb{R}^n)$ for some $0<p<\infty$ \eqref{Boundf}).
Then choosing 
\begin{equation}
\label{muregC}
    \mu^2= \disp \left(\frac{\delta}{C}\right)^{\frac{2}{p+2}},
\end{equation} 
we have 
\begin{equation}
\|  f -f_{\delta ,\mu } \| _{L^2(\mathbb{R}^n)} \le 2 \delta^{\frac{p}{p+2}} C^{\frac{2}{p+2}}\left[M+\frac{1}{2} \max\left\{1;\left(\frac{\delta}{C}\right)^{\frac{2-p}{p+2}}\right\}\right] . 
\end{equation}
\end{thm}
\begin{pf}
From now on we denote $\|   \cdot \|   = \|   \cdot \|  _{L^2(\mathbb{R}^n)}$. \\
Defining $\disp \hat{f}_{\mu}(\bxi) :=  \frac{\Lambda(\bxi)}{1+\mu^2\|\bxi\|^2} \hat{y}(\bxi) $,
 by \eqref{illf}  it follows that  
\begin{equation}
\label{ineqteo2}
\begin{split}
\left| \hat{f}(\bxi )-\hat{f}_{\mu}(\bxi )   \right|  
=&  \left|    \hat{f}(\bxi) \left(1-\frac{1}{1+\|\bxi\|^2\mu^2}\right) \frac{(1+\|\bxi\|^2)^{\frac{p}{2}}}{(1+\|\bxi\|^2)^{\frac{p}{2}}}   \right| \\
   \le & \sup_{\|\bxi\| \in \mathbb{R}} \left| (1+\|\bxi\|^2)^{-\frac{p}{2}}\left(1-\frac{1}{1+\|\bxi\|^2\mu^2}\right)    \right| 		\left|  \hat{f}(\bxi) (1+\|\bxi\|^2)^{\frac{p}{2}} \right|,
\end{split}
\end{equation}
and by the definition of the $H^p(\mathbb{R}^n)$-norm given in \eqref{Boundf}, we have
\begin{equation}
\label{ineq2teo2}
\| \hat{f}(\bxi )-\hat{f}_{\mu}(\bxi )   \|  
 \le  \sup_{\|\bxi\| \in \mathbb{R}} \left| (1+\|\bxi\|^2)^{-\frac{p}{2}}\left(1-\frac{1}{1+\|\bxi\|^2\mu^2}\right)    \right| 		 \|  f \| _{H^p(\mathbb{R}^n)}.
\end{equation}
From the triangle inequality, 
\begin{equation}
\label{dnormas}
 \|   \hat{f} - \hat{f}_{\delta ,\mu }  \|   \leq \|   \hat{f} - \hat{f}_{\mu }  \|    +  \|  \hat{f}_{ \mu }   - \hat{f}_{\delta ,\mu }\|,
\end{equation}
then (\ref{ineq2teo2})-(\ref{dnormas}) and the definition of the regularized source \eqref{ffregtransf}
yield to
\begin{center}
\begin{tabular}{rl}
$\left\|   \hat{f}  -\hat{f}_{\delta ,\mu } \right\| \le$& $\sup\limits_{\|\bxi\| \in \mathbb{R}}  \left|(1+\|\bxi\|^2)^{-\frac{p}{2}}\left(1-\frac{1}{1+\|\bxi\|^2\mu^2}\right) \right| \|  f \| _{H^p(\mathbb{R}^n)} $\\
&\\
&$ +  \sup\limits_{\|\bxi\| \in \mathbb{R}} \left|  \frac{\Lambda (\bxi)}{1+\|\bxi\|^2 \mu^2} \right|
\left\|    \hat{y}-\hat{y}_{\delta } \right\|. $
\end{tabular}
\end{center}
From \cite{YF10},  
$ \disp \sup_{\|\bxi\| \varepsilon \mathbb{R}} \left\vert (1+\|\bxi\|^2)^{-\frac{p}{2}} \left(1-\frac{1}{1+\|\bxi\|^2 \mu^2} \right)\right\vert \le \max \left\{\mu^p,\mu^2\right\},$ this result
together with Lemma \ref{lemma4} and the assumption  $\left\|    \hat{y}-\hat{y}_{\delta }\right\| \leq \delta$ lead to  
\begin{equation*}
\|   \hat{f}  -\hat{f}_{\delta ,\mu }  \| 
\le   \max \left\{ \mu ^{p},\mu^2 \right\} 
       \|  f \| _{H^p(\mathbb{R}^n)} +  \frac{2\delta }{\mu^2} M .
\end{equation*}  
Choosing  $ \mu^2= \disp \left(\frac{\delta}{C}\right)^{\frac{2}{p+2}}$, by Parseval's identity and the linearity of the Fourier transform, 
\begin{equation*}
 \|  f -f_{\delta, \mu} \| =\left\|   \hat{f}  -\hat{f}_{\delta ,\mu } \right\| \le \max \left\{ \delta^{\frac{p}{p+2}} C^{\frac{2}{p+2}};\delta^{\frac{2}{p+2}} C^{\frac{p}{p+2}}\right\} +2M \delta^{\frac{p}{p+2}} C^{\frac{2}{p+2}},
   \end{equation*}
equivalently
\begin{equation}
\label{cota1}
\|  f -f_{\delta ,\mu } \| _{L^2(\mathbb{R}^n)} \le 2 \delta^{\frac{p}{p+2}} C^{\frac{2}{p+2}}\left[M+\frac{1}{2} \max\left\{1;\left(\frac{\delta}{C}\right)^{\frac{2-p}{p+2}}\right\}\right]. 
\end{equation}
\end{pf}

A-priori regularization parameter is usually chosen to be dependent on a-prior bound for the $H^p$ norm of the source and data noise.
For the numerical examples the bound is generally assumed to be 1 \cite{QianFu06,QianFu07,Xiao11,YangFu2012,ZhaoMeng2014} which can lead to erroneous estimates when $\|  f\|  _{H^{p}(\mathbb{R}^{n})}>1$.

Here, we include a rule of choice for the regularization parameter that only depends on data noise. The following theorem is aimed to the estimate error for this case.

\setlength{\leftskip}{0pt}
\begin{thm}
\label{boundestimate}
Consider the inverse problem of determining the source $f(\bx)$ in (\ref{transpeqn})-(\ref{noiselevel}). Let $f_{\delta ,\mu }(\bx)$ be the regularization solution  given in (\ref{ffreg}) and assume that $\|  f \| _{H^p(\mathbb{R}^n)} $ is bounded in $H^p(\mathbb{R}^n)$ for some $0<p<\infty$\eqref{Boundf}.
Then choosing 
\begin{equation}
\label{mureg}
    \mu^2= \disp \delta^{\frac{2}{p+2}},
\end{equation} 
there exists a constant $K$  independent of $\delta$ such that
\begin{equation}
\|  f -f_{\delta ,\mu } \| _{L^2(\mathbb{R}^n)} \le  K \, \max \left\{ \delta^{\frac{2}{p+2}}; \delta^{\frac{p}{p+2}}\right\}. 
\end{equation}
\end{thm}

\begin{pf}
From now on we denote $\|   \cdot \|   = \|   \cdot \|  _{L^2(\mathbb{R}^n)}$. From the proof of Theorem \ref{boundestimate1} we have that
$\|   \hat{f}  -\hat{f}_{\delta ,\mu }  \| 
\le   \max \left\{ \mu ^{p},\mu^2 \right\} 
       \|  f \| _{H^p(\mathbb{R}^n)} +  \frac{2\delta }{\mu^2} M.$
Choosing  $ \mu^2= \disp \delta^{\frac{2}{p+2}}$, by Parseval's identity and the linearity of the Fourier transform, 
\begin{equation}
\label{cota}
 \|  f -f_{\delta, \mu} \| =\left\|   \hat{f}  -\hat{f}_{\delta ,\mu } \right\| \le K  \max \left\{ \delta^{\frac{2}{p+2}}, \delta^{\frac{p}{p+2}}\right\}.
 \end{equation}
 where $K= C + 2M= C + 2 \max \left\{\frac{1}{t_0} + \frac{\sqrt{n} \|\bbeta\|_\infty}{2\nu t_0} ;\nu +\alpha^2+\frac{\sqrt{n} \|\bbeta\|_\infty}{2} \right\}$ and $C$ is the bound for the $H^p(\mathbb{R}^n)$-norm of $f$, i.e., $\|  f\|  _{H^{p}(\mathbb{R}^n)} \le C$.
\end{pf}

\begin{rem}
Note $p = \infty$ is excluded since and the error bound in that case the error bound is $K$.

A particular case of mathematical interest is for $p = 2$, where we obtain
$$ \|  f -f_{\delta, \mu} \|  \le K  \sqrt{\delta}.$$

\end{rem}

\begin{rem}
If a bound ${\delta_M}>1$  is allowed for noise in measurements,  in order to keep $0<\mu<1$ one can take 
\begin{equation}\label{PRII}
\mu^2= \disp \left(\frac{\delta}{\delta_M}\right)^{\frac{2}{p+2}}.
\end{equation}

In that case, for $K= C + 2 \, \disp {\delta_M}^{\frac{4}{p+2}}\max \left\{\frac{1}{t_0} + \frac{\sqrt{n} \|\bbeta\|_\infty}{2\nu t_0} ;\nu +\alpha^2+\frac{\sqrt{n} \|\bbeta\|_\infty}{2} \right\}$ we have
 $\|  f -f_{\delta, \mu} \| =\left\|   \hat{f}  -\hat{f}_{\delta ,\mu } \right\| \le K  \max \left\{ \left(\frac{\delta}{\delta_M}\right)^{\frac{2}{p+2}}; \left(\frac{\delta}{\delta_M}\right)^{\frac{p}{p+2}}\right\}.$
\end{rem}

\section{Numerical examples}\label{s:4} \vspace{-4pt}

In this section  we consider few examples with a source $f\in \mathbb{R}^n,  n=1,2,3$  to illustrate the performance of the regularization operator. 
For each of them we have chosen different values for the parameters $ \alpha^2,\bbeta,\nu$, $t_0$ and a set of standard deviation values $\{\epsilon_1,...,\epsilon_k\}$ for the data noise.
The space is uniformly discretized and a data set $\{y_{\delta_1},...,y_{\delta_N}\}$ is obtained by evaluating the solution $u(\bx,t)$ at a fixed time instant $t_0$ and adding noise, that is, 
\begin{equation}
\label{noisydata}
y_{\delta_i} = y(\bx_i) + \epsilon_i , \quad i=1,...,N,\,\,  \bx_i \in {\cal G},
\end{equation}
where ${\cal G}$ is the uniform grid defined on $\mathbb{R}^n$ and $\epsilon_i, i=1,...,N$ are realizations of the normally distributed random variable ${\eta}$ with mean 0 and standard deviation $\epsilon$. 

By denoting $y_i=y(\bx_i),  i=1,.., N$, since the noise level  $\delta$ satisfies \eqref{noiselevel}, the error 
\begin{equation}
\label{discrerror}
 \by -\by_{\delta} =(y_1 - y_{\delta_1},..., y_N-y_{\delta_N})=(\epsilon_1,...,\epsilon_N),
 \end{equation}
 is numerically integrated using the Simpson's method to obtained an approximated value for $\delta=\delta(\epsilon)$. Then $\delta_M$ in \eqref{noiselevel}  is chosen to be an upper bound for $\delta$. In practice, $\delta_M$ can be estimated based on the measuring instruments used in data collection. 

Afterwards,$\{ \hat{y}_{\delta_1},...,\hat{y}_{\delta_N} \}$ is calculated by means of the FFT (Fast Fourier Transform) and the regularized solution $f_{\delta ,\mu}$ is calculated as defined in \eqref{ffreg}  using the inverse FFT \cite{Elden00} where the regularization parameter $\mu$ is chosen according to \eqref{PRII}.

The results of the estimated sources (non-regularized and the regularized one) are plotted. A table of errors is included for each example which shows the absolute and relative errors. For comparison purposes, the errors for all the examples are calculated taking \eqref{discrerror} $\{\epsilon_1,...,\epsilon_5\}= \{0.01, 0.03, 0.05, 0.08, 0.1\}$. 

\begin{rem}
Notice that in some examples a good estimation is obtained even when the parameter value $p$  does not correspond to the $H^p$ space where the source belongs.
\end{rem}

\subsection{Examples 1D}

Two examples are considered for the one-dimensional inverse source problem defined in \eqref{transpeqn}. 

\begin{exmp}
\label{example1}
For this example the source $f$ is defined by 
$$f(x)=\begin{cases} 
-1,\qquad & -20 \leq x<-10, \\
1,\qquad & -10 \leq x<0, \\
-1, \qquad & 0 \leq x<10, \\
1, \qquad & 10 \leq x \leq 20, \\
0, \qquad & \text {in another case}.
\end{cases}$$
With modeling parameter values  $\alpha^2=2.10^{-5}, \bbeta= 1.10^{-5}, \nu= 1, t_0=5, p=1$ and $\epsilon \in \{0.2 , 0.15, 0.1, 0.05\}$. 
\end{exmp}

\begin{figure}[h!]
\centering {
\begin{tabular}{ccc}
{\includegraphics[ scale=0.48]{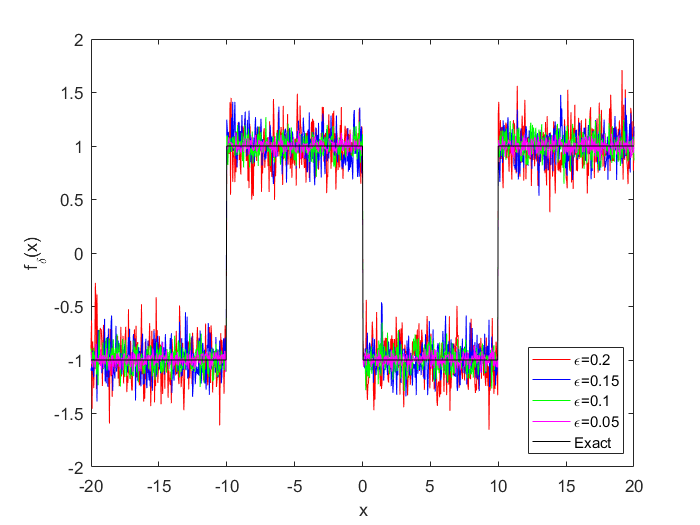}}&
{\includegraphics[ scale=0.48]{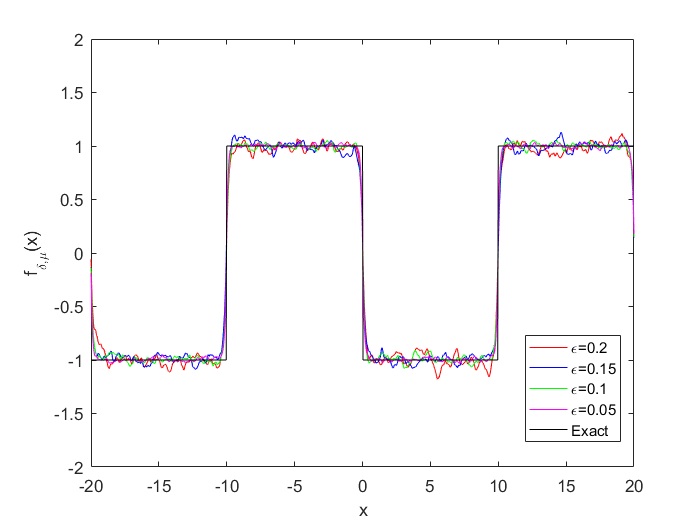}}\\
(a) & (b)\\
\end{tabular}
\caption{Sources for example \ref{example1}: unregularized (a) and regularized (b) assuming $\alpha^2=2.10^{-5}, \bbeta= 1.10^{-5}, \nu= 1, t_0=5, p=1$ and $\epsilon \in \{0.2 , 0.15, 0.1, 0.05\}$.}\label{CurvaCuadrada}
}
\end{figure}

\begin{table}[h!]
\begin{center}
\caption{Example \ref{example1}:  Errors assuming $\alpha^2=2.10^{-5}, \bbeta= 1.10^{-5}, \nu= 1, t_0=1, p=1$.}

{\begin{tabular}{cc}\hline
Absolute errors & Relative errors\\
{\begin{tabular}{lccc} \hline
$\epsilon$ & $\|  f-f_{\delta} \|$ & $\|  f-f_{\delta,\mu} \|$ \\ \hline
0.01  & 0.4610  & 0.0998     \\
0.03  & 0.5750  & 0.2985      \\
0.05  & 0.6430  &	0.4853    \\
0.08  & 0.8268  & 0.7000     \\
0.1   & 0.9512  & 0.7445    \\\hline
\end{tabular}}&
{\begin{tabular}{cc} \hline
$\|  f-f_{\delta} \|/\|  f \|$ & $\|  f-f_{\delta,\mu} \|/\|  f \|$  \\ \hline
0.0729  & 0.0158\\
0.0909  & 0.0472\\
0.1017	& 0.0767\\
0.1307  & 0.1107\\
0.1504  & 0.1177 \\ \hline
\end{tabular}}
\end{tabular}}
\label{tableej1}
\end{center}
\end{table}

\bigskip
\medskip


\begin{exmp}
\label{example2}
For this example the source $f$ is defined by 
$$f(x)=\begin{cases} 
x+1, \qquad & -1 \leq x < 0, \\
-x+1, \qquad & 0 \leq x \leq 1, \\
0, \qquad & \text {in another case}.
\end{cases}$$
With modeling parameter values  $\alpha^2=2, \bbeta= 0, \nu= 1, t_0=0.2, p=2$ and $\epsilon \in \{0.004 , 0.003, 0.002, 0.001\}$. 

\end{exmp}

\begin{figure}[h!]
\centering {
\begin{tabular}{ccc}
{\includegraphics[ scale=0.48]{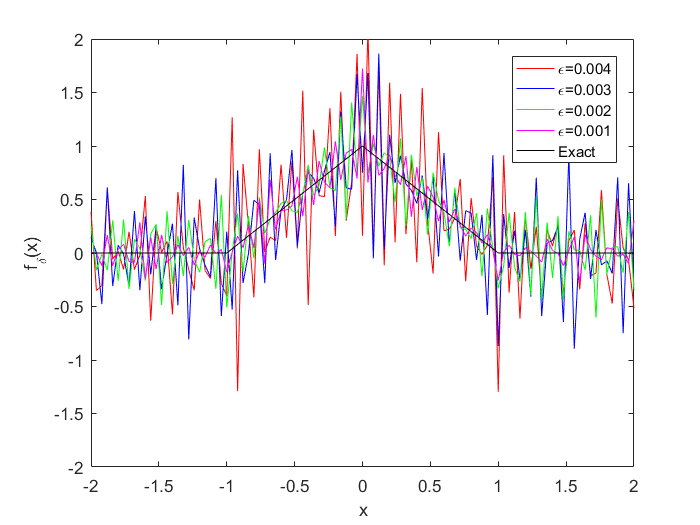}}&
{\includegraphics[ scale=0.48]{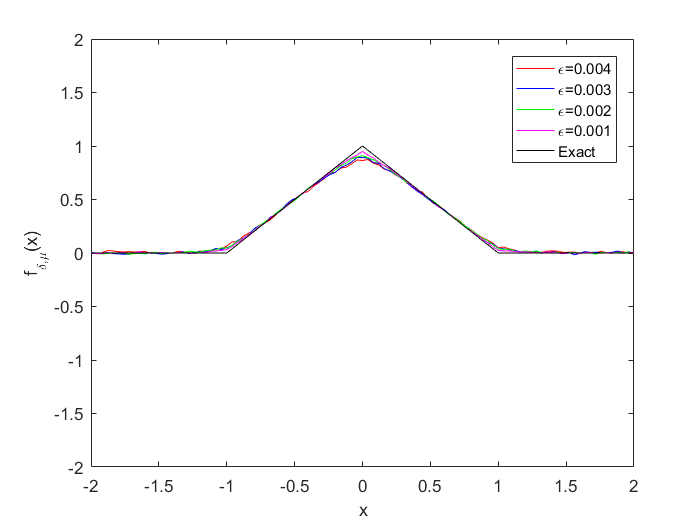}}\\
(a) & (b)\\
\end{tabular}
\caption{Sources for example \ref{example2}: unregularized (a) and regularized (b) assuming $\alpha^2=2, \bbeta= 0, \nu= 1, t_0=0.2, p=2$ and $\epsilon \in \{0.004 , 0.003, 0.002, 0.001\}$.}\label{CurvaTriangular}
}
\end{figure}

\begin{table}[h!]
\begin{center}
\caption{Example \ref{example2}: Errors assuming $\alpha^2=2, \bbeta= 0, \nu= 1, t_0=1, p=1$.}
{\begin{tabular}{cc} \hline
Absolute errors & Relative errors \\
{\begin{tabular}{lccc} \hline
$\epsilon$ & $\|  f-f_{\delta} \|$ & $\|  f-f_{\delta,\mu} \|$ \\ \hline
0.01  & 2.5927   & 0.1698   \\
0.03  & 9.4906   & 0.2559    \\
0.05  & 12.9565  &	0.2878 \\
0.08  & 23.9564  & 0.3579   \\
0.1   & 25.1220  & 0.3713   \\\hline
\end{tabular}}&
{\begin{tabular}{cc}\hline
 $\|  f-f_{\delta} \|/\|  f \|$ & $\|  f-f_{\delta,\mu} \|/\|  f \|$  \\ \hline
 3.1728   & 0.2077  \\
 11.6143  & 0.3131  \\
 15.8558	 & 0.3522	 \\
 29.7091  & 0.4380  \\
 30.7802  & 0.4512  \\\hline
\end{tabular}}
\end{tabular}}
\label{tableej2}
\end{center}
\end{table}

\subsection{Examples 2D}

Two examples are considered for the two-dimensional inverse source problem defined in \eqref{transpeqn}. 


\begin{exmp}
\label{example3}

For this example the source $f$ is defined by 
$$f(x)=\begin{cases} 
\ cos\left(\frac{x_1}{20}\right) \ cos\left(\frac{x_2}{20}\right), \qquad & -40 \leq x_1,x_2 \leq 40, \\
0, \qquad & \text {in another case}.
\end{cases}$$
With modeling parameter values  $\alpha^2=0.2, \bbeta= (0,0), \nu= 0.999, t_0=1, p=1$ and $\epsilon=0.025$. 

\end{exmp}

\begin{figure}[h!]
\centering {
\begin{tabular}{ccc}
{\includegraphics[ scale=0.30]{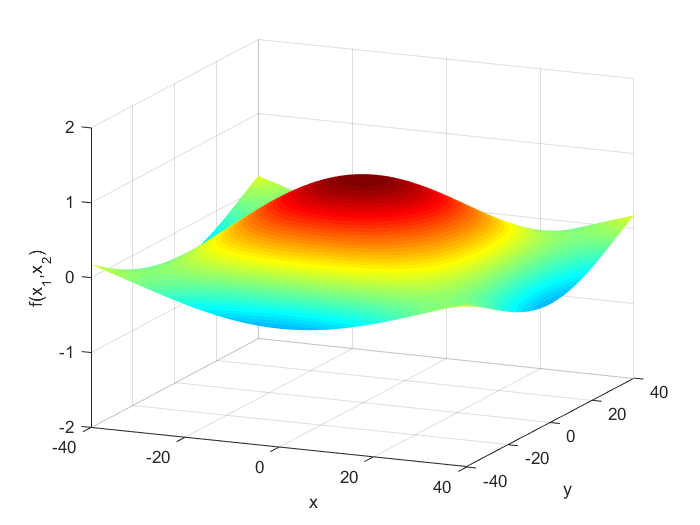}}&
{\includegraphics[ scale=0.30]{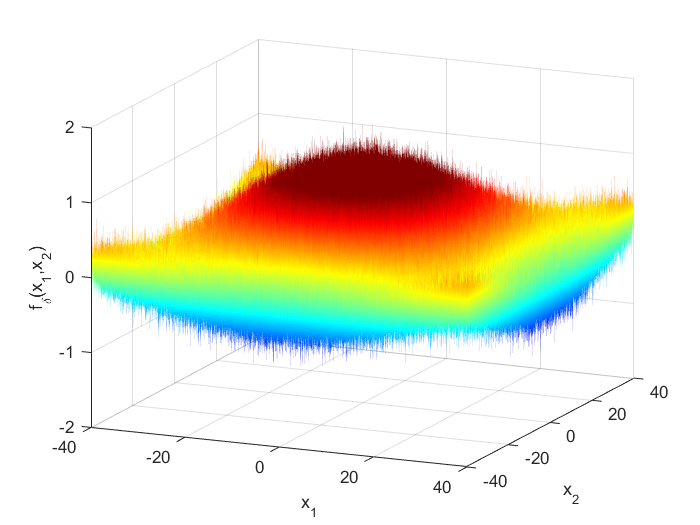}}&
{\includegraphics[ scale=0.30]{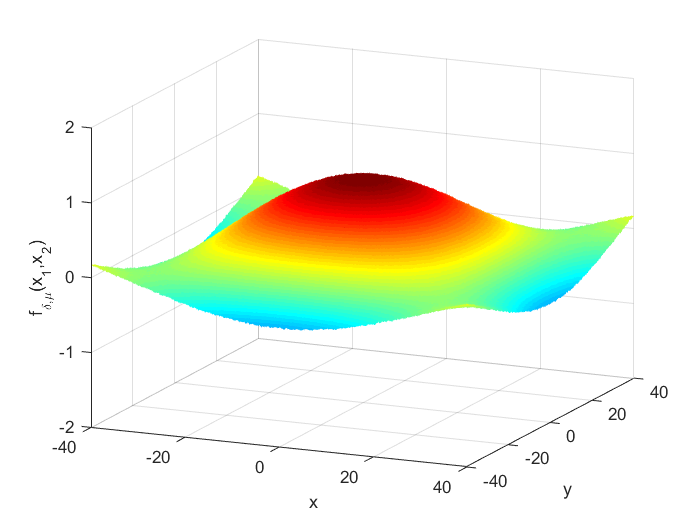}}\\
{\includegraphics[ scale=0.30]{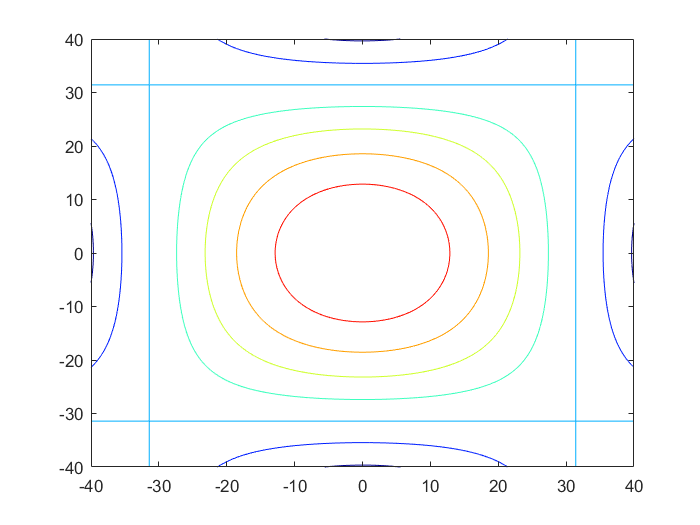}}&
{\includegraphics[ scale=0.30]{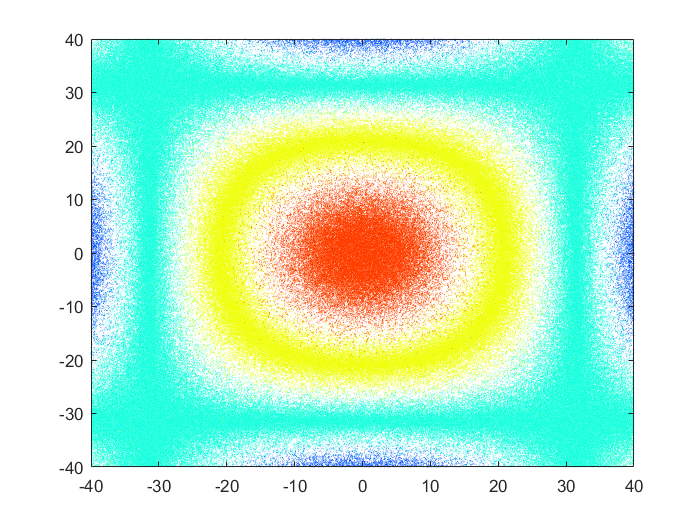}}&
{\includegraphics[ scale=0.30]{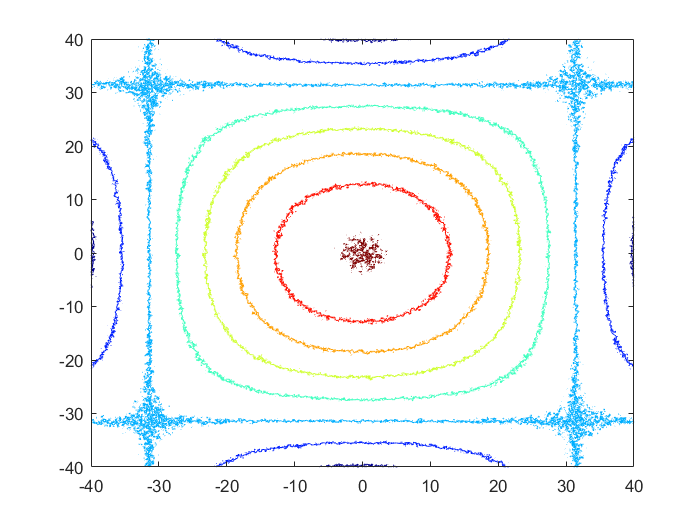}}\\
(a) & (b) & (c)\\
\end{tabular}
\caption{Sources for example \ref{example3}: original (a), unregularized (b) and regularized (c) assuming $\alpha^2=0.2, \bbeta= (0,0), \nu= 0.999, t_0=1, p=1$ and $\epsilon=0.025$.}\label{SuperficieSenosoidal}
}
\end{figure}

\begin{table}[ht!]
\begin{center}
\caption{Example \ref{example3}: Errors assuming  $\alpha^2=0.2, \bbeta= (0,0), \nu= 0.999, t_0=1, p=1$.}
{\begin{tabular}{cc} \hline
Aboslute errors & Relative errors \\
{\begin{tabular}{lccc} \hline
$\epsilon$ & $\|  f-f_{\delta} \|$ & $\|  f-f_{\delta,\mu} \|$   \\ \hline
0.01  & 5.6563    & 0.2599      \\
0.03  & 16.9395   & 0.7778   \\
0.05  & 28.2357   &	1.3021   	\\
0.08  & 45.1974   & 2.0774    \\
0.1   & 56.3882   & 5.5962    \\\hline
\end{tabular}}&
{\begin{tabular}{cc} \hline
$\|  f-f_{\delta} \|/\|  f \|$ & $\|  f-f_{\delta,\mu} \|/\|  f \|$  \\ \hline
 0.1744   & 0.0008  \\
0.5221   & 0.0240  \\
0.8706	 & 0.0401	 \\
1.3936   & 0.0641  \\
1.7387   & 0.0800  \\\hline
\end{tabular}}
\end{tabular}}
\label{tableej3}
\end{center}
\end{table}

\medskip
\bigskip
\bigskip


\begin{exmp}
\label{example4}
For this example the source $f$ is defined by 
$$f(x_1,x_2)=\begin{cases} 
10+x_1-x_2, \qquad & -10 \leq x_1 \leq 0, \,\,\,\,\,\, 0 \leq x_2 \leq 10+x_1 , \\
10+x_1+x_2, \qquad & -10 \leq x_1 \leq 0, \,\,\,\,\,\, -10-x_1 \leq x_2 \leq 0 , \\
10-x_1-x_2, \qquad & 0 \leq x_1 \leq 10,  \,\,\,\,\,\,  0 \leq x_2 \leq 10-x_1 , \\
10-x_1+x_2, \qquad & 10 \leq x_1 \leq 10, \,\,\, \,\,\, -10+x_1 \leq x_2 \leq 0 , \\
0, \qquad & \text {in another case}.
\end{cases}$$

With modeling parameter values  $\alpha^2=1, \bbeta= (0,0), \nu= 1, t_0=0.4, p=0.6$ and $\epsilon=0.05$.

\end{exmp}

\begin{figure}[h!]
\centering {
\begin{tabular}{ccc}
{\includegraphics[ scale=0.30]{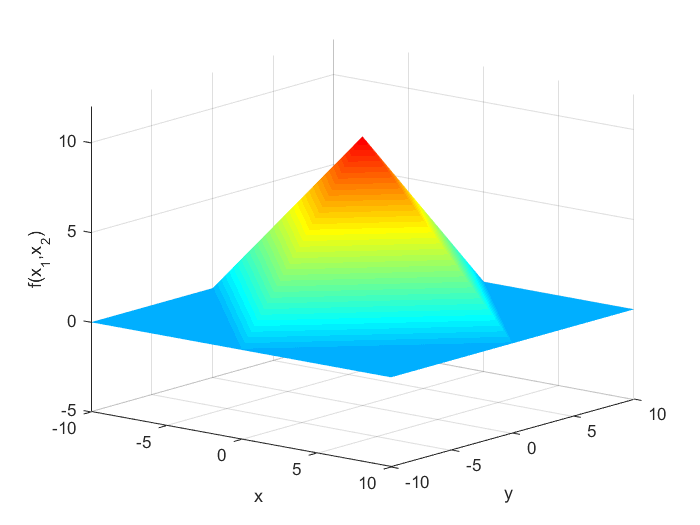}}&
{\includegraphics[ scale=0.30]{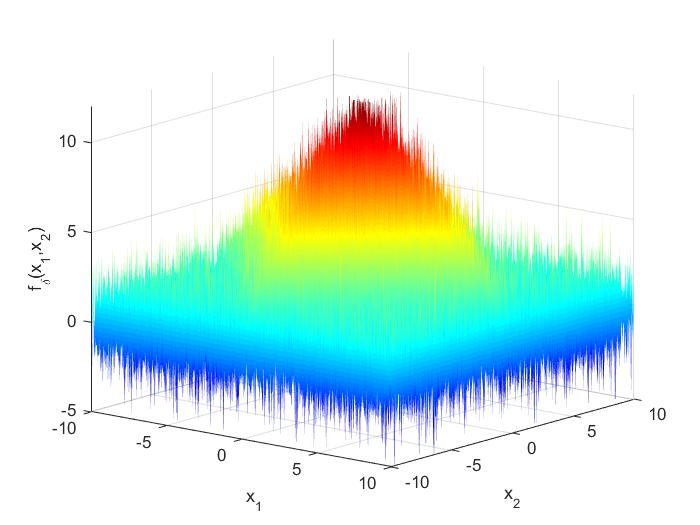}}&
{\includegraphics[ scale=0.30]{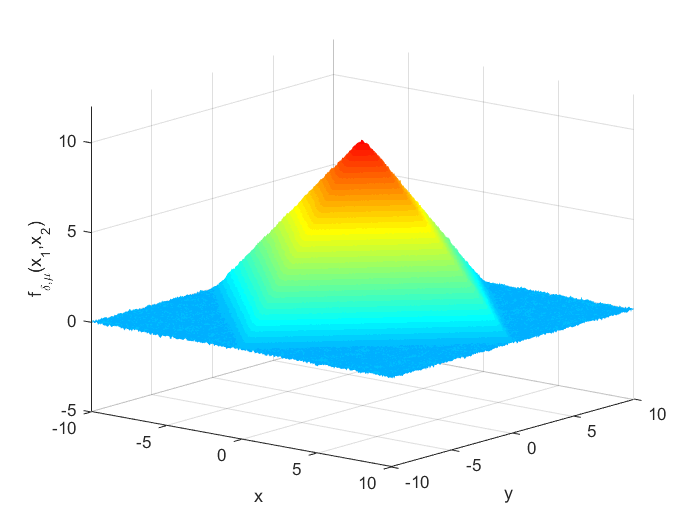}}\\
{\includegraphics[ scale=0.30]{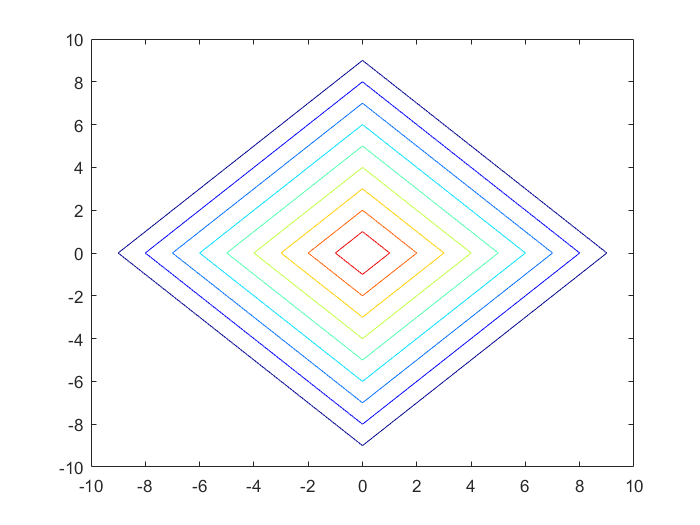}}&
{\includegraphics[ scale=0.30]{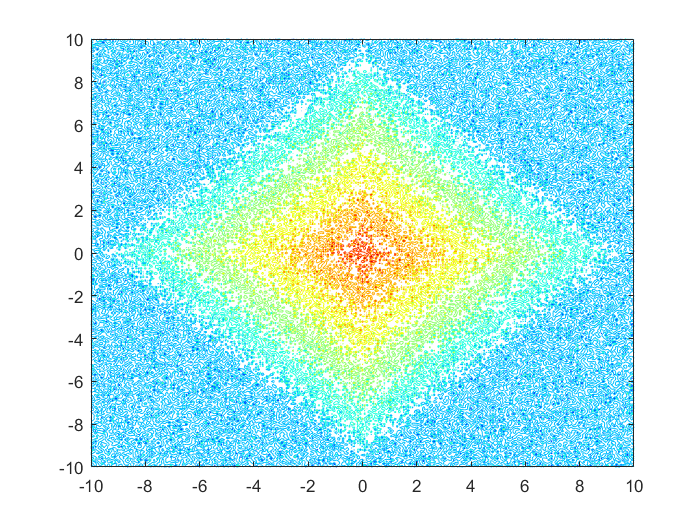}}&
{\includegraphics[ scale=0.30]{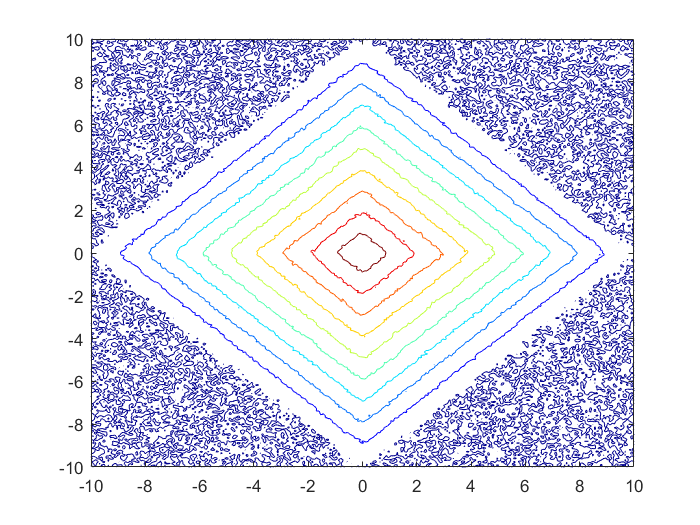}}\\
(a) & (b) & (c)\\
\end{tabular}
\caption{Sources for example \ref{example4}: original (a), unregularized (b) and regularized (c) assuming $\alpha^2=1, \bbeta= (0,0), \nu= 1, t_0=0.4, p=0.6$ and $\epsilon=0.05$.}\label{SuperficiePrisma}
}
\end{figure}

\begin{table}[h!]
\begin{center}
\caption{Example \ref{example4}: Errors assuming $\alpha^2=1, \bbeta= (0,0), \nu= 1, t_0=1, p=1$.}
{\begin{tabular}{cc} \hline
Aboslute errors & Relative errors \\
{\begin{tabular}{lccc} \hline
$\epsilon$ & $\|  f-f_{\delta} \|$ & $\|  f-f_{\delta,\mu} \|$   \\ \hline
0.01  & 6.5149    & 0.3527      \\
0.03  & 19.2785   & 0.6748   \\
0.05  & 31.7283   &	1.0502   	\\
0.08  & 50.7739   & 1.6429    \\
0.1   & 63.9207   & 2.0451    \\\hline
\end{tabular}}&
{\begin{tabular}{cc} \hline
$\|  f-f_{\delta} \|/\|  f \|$ & $\|  f-f_{\delta,\mu} \|/\|  f \|$  \\ \hline
 0.1128   & 0.0061  \\
 0.3339   & 0.0117  \\
0.5496	  & 0.0182	 \\
0.8796   & 0.0285  \\
 1.1071   & 0.0354  \\\hline
\end{tabular}}
\end{tabular}}
\label{tableej4}
\end{center}
\end{table}

\bigskip

\subsection{Examples 3D}

One example are considered for the three-dimensional inverse source problem defined in \eqref{transpeqn}. 


\begin{exmp}
\label{example5}

For this example the source $f$ is defined by 
$$f(x)=\begin{cases} 
\ sin\left(\frac{x_1+x_2+x_3}{20}\right), \qquad & -2\pi \leq x_1,x_2,x_3 \leq 2\pi, \\
0, \qquad & \text {in another case}.
\end{cases}$$
With modeling parameter values  $\alpha^2=0.4, \bbeta= (1,-0.5,-0.5), \nu= 0.997, t_0=3, p=3$ and $\epsilon=0.035$. For this case, we consider
either $x=0$, $y=0$ or $z=0$ to plot the resulting  estimated sources, the non-regularized and the regularized one. 
 
\end{exmp}

\begin{figure}[h!]
\centering {
\begin{tabular}{ccc}
{\includegraphics[ scale=0.30]{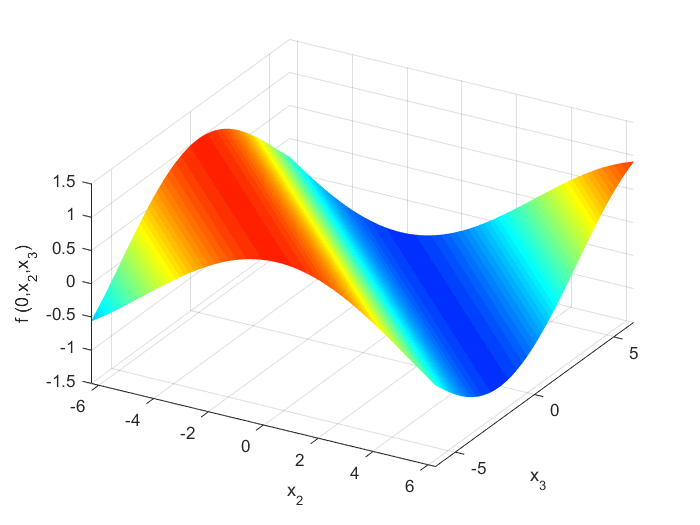}}&
{\includegraphics[ scale=0.30]{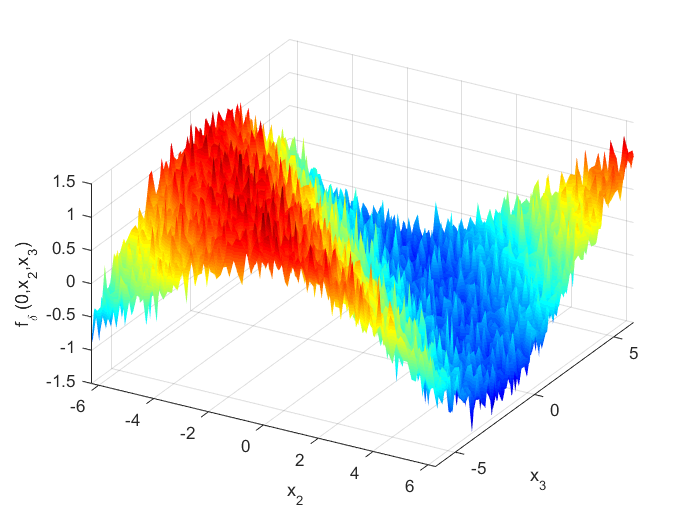}}&
{\includegraphics[ scale=0.30]{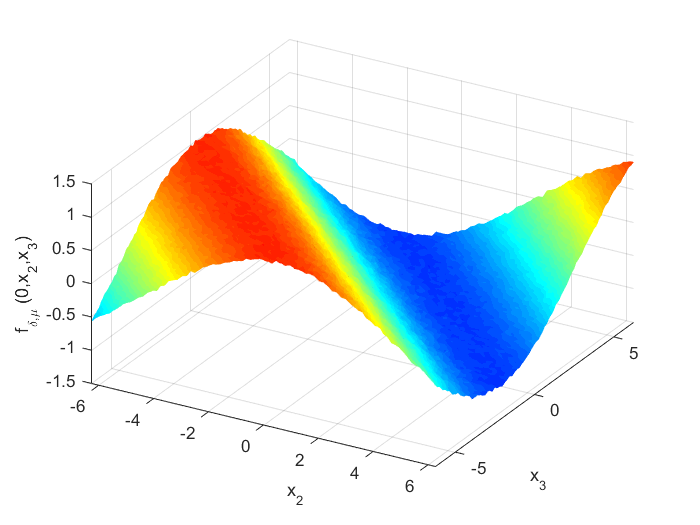}}\\
{\includegraphics[ scale=0.30]{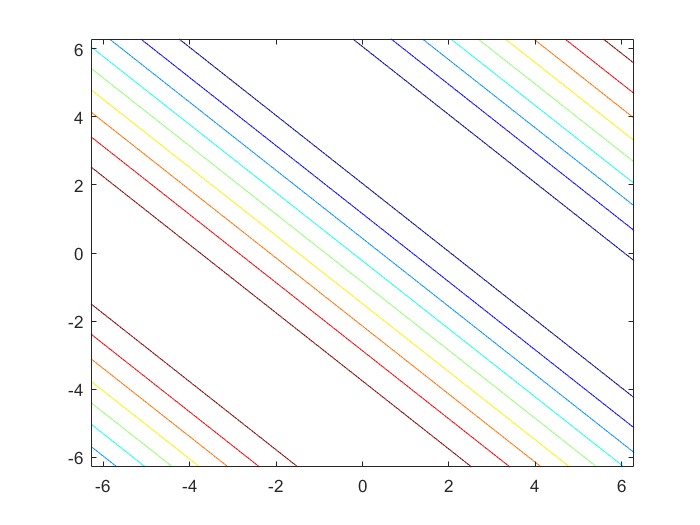}}&
{\includegraphics[ scale=0.30]{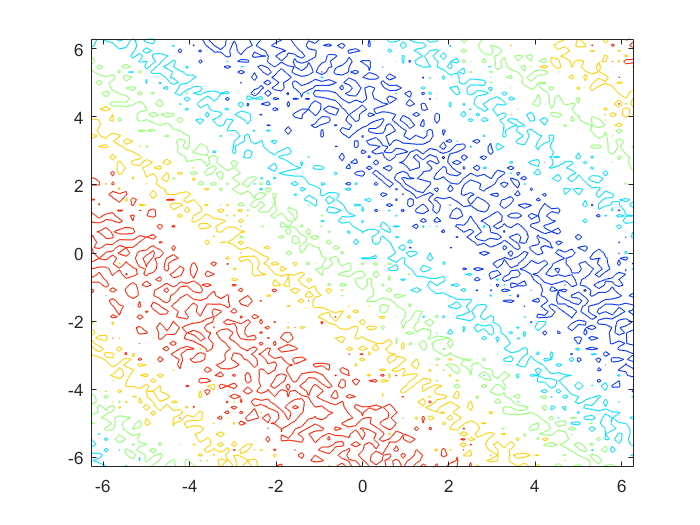}}&
{\includegraphics[ scale=0.30]{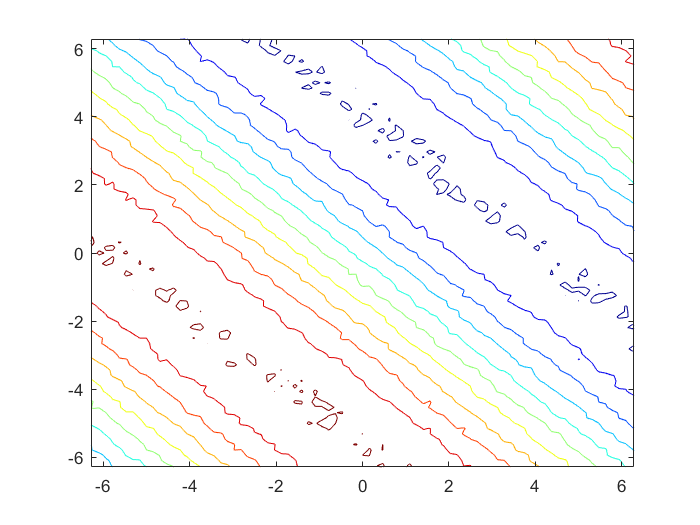}}\\
(a) & (b) & (c)\\
\end{tabular}
\caption{Sources for example \ref{example5}: original (a), unregularized (b) and regularized (c) for $x=0$ assuming $\alpha^2=0.4, \bbeta= (1,-0.5,-0.5), \nu= 0.997, t_0=3, p=3$ and $\epsilon=0.035$.}\label{Volumen_Senosoidal_x_0}
}
\end{figure}

\begin{figure}[h!]
\centering {
\begin{tabular}{ccc}
{\includegraphics[ scale=0.30]{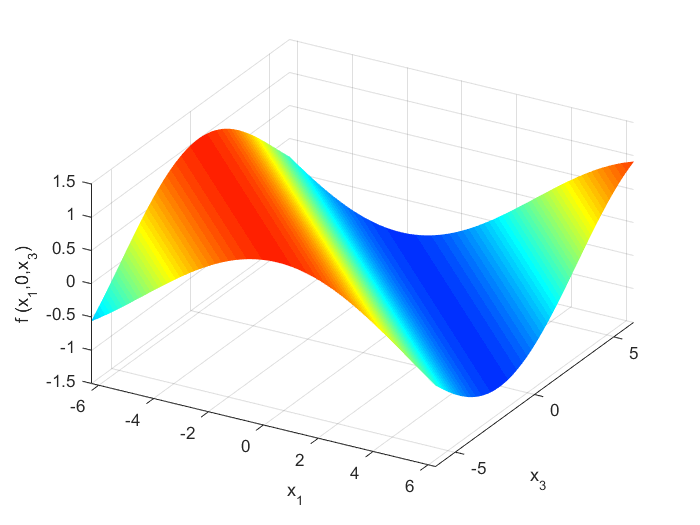}}&
{\includegraphics[ scale=0.30]{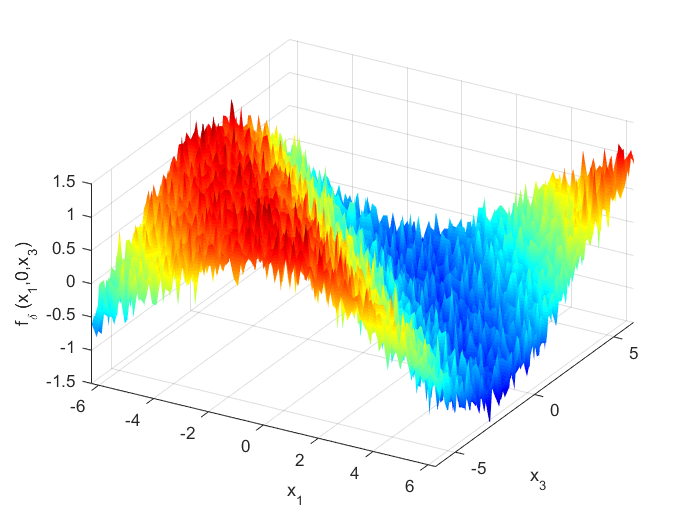}}&
{\includegraphics[ scale=0.30]{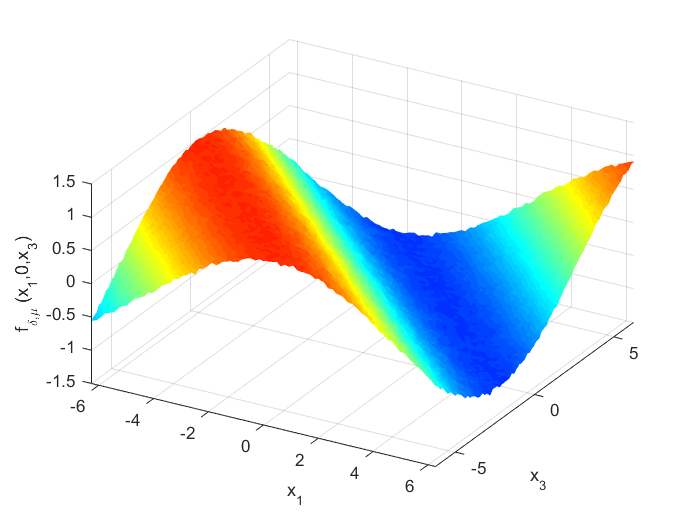}}\\
{\includegraphics[ scale=0.30]{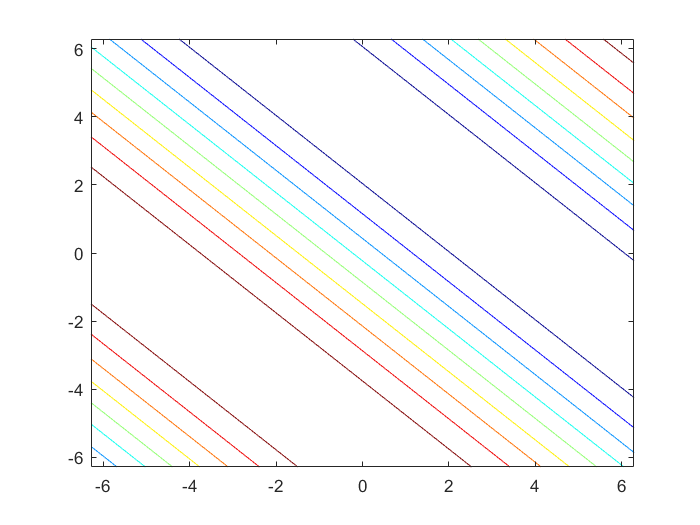}}&
{\includegraphics[ scale=0.30]{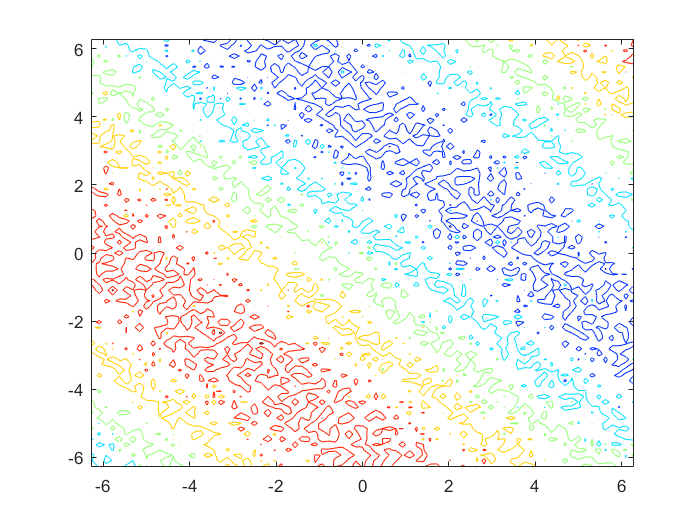}}&
{\includegraphics[ scale=0.30]{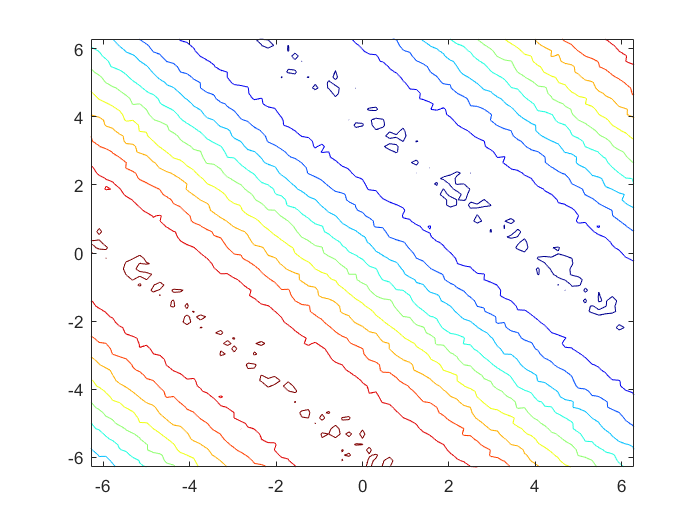}}\\
(a) & (b) & (c)\\
\end{tabular}
\caption{Sources for example \ref{example5}: original (a), unregularized (b) and regularized (c) for $y=0$ assuming $\alpha^2=0.4, \bbeta= (1,-0.5,-0.5), \nu= 0.997, t_0=3, p=3$ and $\epsilon=0.035$.}\label{Volumen_Senosoidal_y_0}
}
\end{figure}

\begin{figure}[h!]
\centering {
\begin{tabular}{ccc}
{\includegraphics[ scale=0.30]{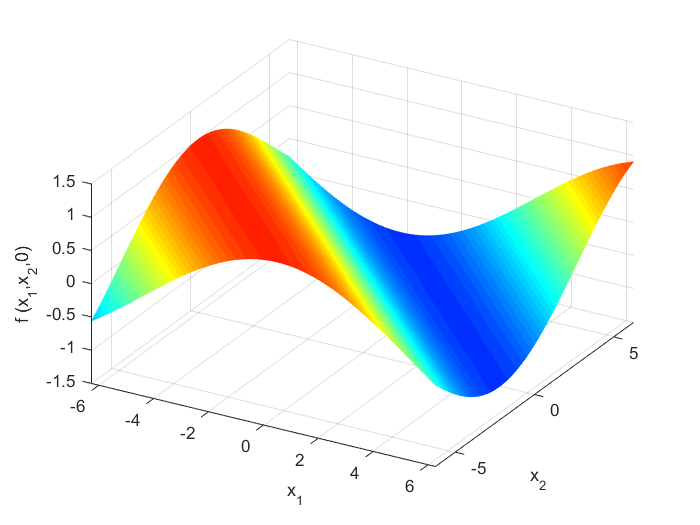}}&
{\includegraphics[ scale=0.30]{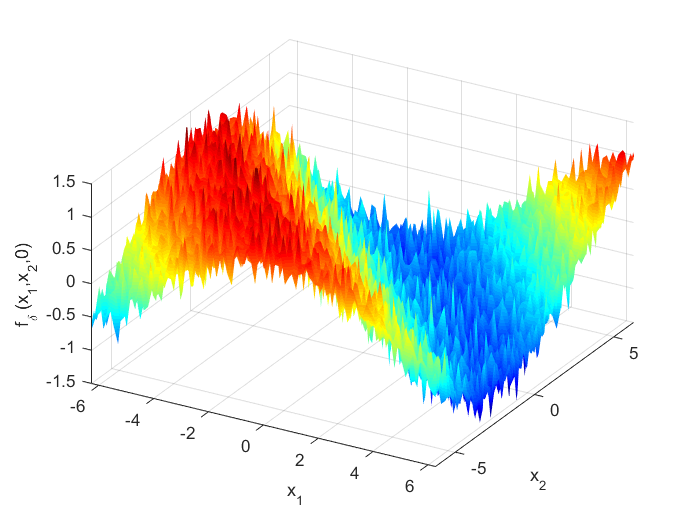}}&
{\includegraphics[ scale=0.30]{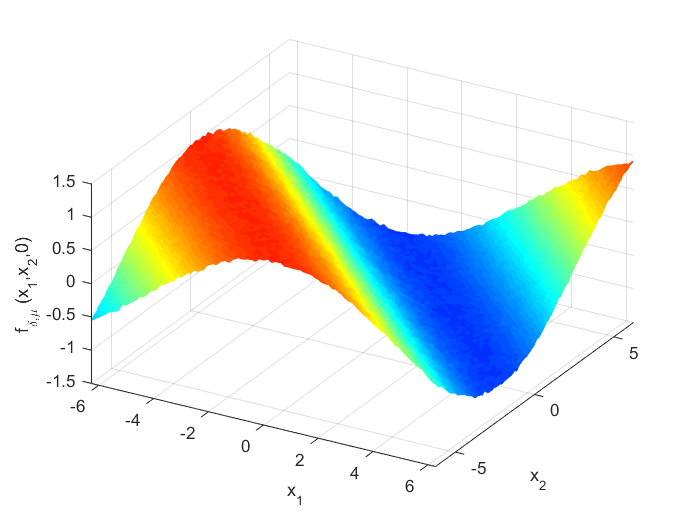}}\\
{\includegraphics[ scale=0.30]{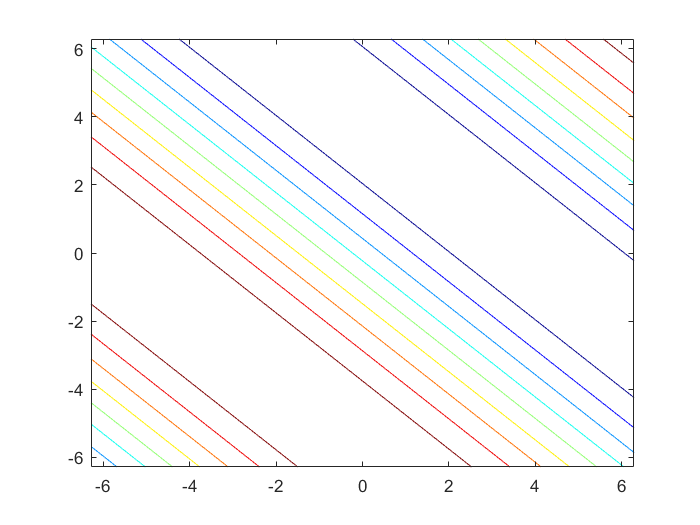}}&
{\includegraphics[ scale=0.30]{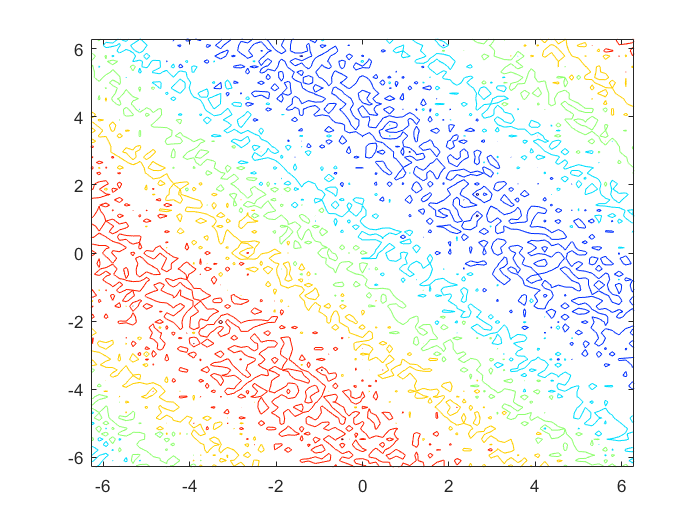}}&
{\includegraphics[ scale=0.30]{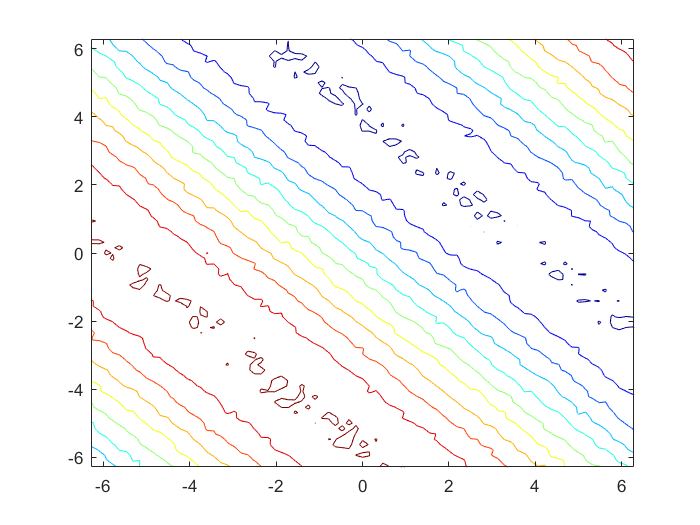}}\\
(a) & (b) & (c)\\
\end{tabular}
\caption{Sources for example \ref{example5}: original (a), unregularized (b) and regularized (c) for $z=0$ assuming $\alpha^2=0.4, \bbeta= (1,-0.5,-0.5), \nu= 0.997, t_0=3, p=3$ and $\epsilon=0.035$.}\label{Volumen_Senosoidal_z_0}
}
\end{figure}

\begin{table}[h!]
\begin{center}
\caption{Example \ref{example5}: Errors assuming $\alpha^2=0.4, \bbeta= (1,-0.5,-0.5), \nu= 0.997, t_0=1, p=1$.}
{\begin{tabular}{cc} \hline
Aboslute errors & Relative errors \\
{\begin{tabular}{lccc} \hline
$\epsilon$ & $\|  f-f_{\delta} \|$ & $\|  f-f_{\delta,\mu} \|$  \\ \hline
0.01  & 2.7377    & 0.3890     \\
0.03  & 8.2057    & 0.7947     \\
0.05  & 13.6597   &	1.2621  \\
0.08  & 21.9170   & 1.9791   \\
0.1   & 27.3875   & 2.4679   \\\hline
\end{tabular}}&
{\begin{tabular}{cc} \hline
$\|  f-f_{\delta} \|/\|  f \|$ & $\|  f-f_{\delta,\mu} \|/\|  f \|$  \\ \hline
 0.0842   & 0.0120  \\
 0.2522   & 0.0244  \\
0.4149	  & 0.0388	 \\
0.6737    & 0.0608  \\
 0.8419   & 0.0759  \\\hline
\end{tabular}}
\end{tabular}}
\label{tableej6}
\end{center}
\end{table}

\section{Conclusion}
\label{S3} \vspace{-4pt}

This work focus on the problem of the inverse source for full parabolic equations in $\mathbb{R}^n.$ A family of regularization operators is defined in order to deal with the ill-posedness the problem. It was designed to compensate the instability factor in the inverse operator. A rule of choice for the regularization parameters is also included which is based on the noise level assumed in the data set and the smoothness of the source to be identified. We demonstrate that for the parameter choice rule proposed here, the method is stable and a bound of H\"older type is obtained for the regularization error.

The numerical examples show good estimates for the different n-variables sources,  $n=1,2,3$, determined from data with different noise levels. Moreover, the sources used for the numerical examples belong to different Hilbert spaces and all of them show a good performance of the regularization approach adopted.

\vspace{10pt} \noindent
{\bf Acknowledgements:}  I thank Diana Rubio for her valuable and selfless collaboration.

\vspace{10pt} \noindent
{\bf Funding:}  This work was partially supported by SOARD/AFOSR through grant FA9550-18-1-0523,  CONICET graduate scholarship and UNGS research assistantship.

\end{document}